\documentclass[12pt]{amsart} 
\usepackage{geometry}
\usepackage{graphicx} 
\usepackage{amssymb}
\usepackage{epstopdf}
\usepackage{amsmath,amscd}
\usepackage{amsthm} 
\usepackage{url,verbatim}
\usepackage{mathtools}

\RequirePackage[colorlinks,citecolor=blue,urlcolor=blue]{hyperref}
\usepackage{breakurl} 
\theoremstyle{plain}
\newtheorem{theorem}{Theorem}
 
\newtheorem{lemma}[theorem]{Lemma}

\newcounter{mycount}

\numberwithin{equation}{section} 
\numberwithin{theorem}{section}
\numberwithin{figure}{section}

\newcommand\RR{{\mathbb R}} 
 
\newcommand\qq{\qquad} 
\newcommand\q{\quad}

\newcommand\om{\omega} 
\newcommand\Om{\Omega}

\newcommand\NN{{\mathbb N}} 
\newcommand\ZZ{{\mathbb Z}}

\newcommand\eps{\epsilon} 
 
\newcommand\resp{respectively}
 
\newcommand\lra{\leftrightarrow}
\newcommand\oo{\infty}

\newcommand\pc{p_{\mathrm{c}}}
\newcommand\wpc{\widehat p_{\mathrm{c}}}

\newcommand\pic{\pi_{\mathrm{c}}}
\newcommand\rhoc{\rho_{\mathrm{c}}}

\newcommand\LR{\mathrm{LR}}
\newcommand\pd{\partial}

\title
{The work of Lucio Russo on percolation}

\author{Geoffrey R.\ Grimmett} \address{Statistical Laboratory, Centre for
Mathematical Sciences, Cambridge University, Wilberforce Road, Cambridge CB3
0WB, UK} 
\email{g.r.grimmett@statslab.cam.ac.uk}
\urladdr{\url{http://www.statslab.cam.ac.uk/~grg/}}

\begin{document}

\begin{abstract} 
The contributions of Lucio Russo to the mathematics of percolation and disordered systems are
outlined. The context of his work is explained, and its ongoing impact
on current work is described and amplified.
\end{abstract}

\date{24 March 2016} 

\dedicatory{Dedicated in friendship to Lucio Russo}
\maketitle

\section{A personal appreciation}\label{sec:pers}

Prior to his mid-career move to the history of science in the early 1990s, Lucio Russo enjoyed a
very successful and influential career in the theory of probability and disordered systems,
in particular of percolation and the Ising model. His ideas have shaped these 
significant fields of science, and his name will always be associated with
a number of fundamental techniques of enduring importance.

The author of this memoir is proud to have known Lucio in those days, and to have profited
from his work, ideas, and company. He hopes that this brief account of some of Lucio's results
will stand as testament to the beauty and impact of his ideas.

\section{Scientific summary}\label{sec:sci}

Lucio Russo has  worked principally on the mathematics of percolation, that is, of the existence (or not)
of infinite connected clusters within a disordered spatial network. The principal
model in this field is the so called
\emph{percolation model}, introduced to mathematicians by Broadbent and Hammersley
in 1957, \cite{BH57}. Consider, for definiteness, the hypercubic lattice $\ZZ^d$ with $d \ge 2$, and 
let $p\in [0,1]$. We
declare each edge to be \emph{open} with probability
$p$ and \emph{closed} otherwise, and different edges receive independent states.
The main questions are centred around the existence (or not) of an infinite open
component in $\ZZ^d$. It turns out that there exists a critical probability $\pc=\pc(d)$
such that no infinite open cluster exists when $p<\pc$, and there exists a unique such cluster when $p>\pc$.
(It is not still known which of these two occurs when $p=\pc$ for general $d$, specifically
when $3\le d \le 10$. See \cite{FvdH}.)

Let $C$ be the open cluster of $\ZZ^d$ containing the origin.
Two  functions that play important roles in the theory are the \emph{percolation probability}
$\theta$ and the \emph{mean cluster size} $\chi$ given by
$$
\theta(p)=P_p(|C|=\oo),\qq \chi(p)=E_p|C|,
$$
where $P_p$ and $E_p$ are the appropriate product measure and expectation.
The above model is the  \emph{bond} percolation
model; the \emph{site} percolation model is defined similarly, with sites being open/closed.
A fairly recent account of percolation may be found in \cite{G99}. 

The question was raised  in 1960 (by Harris, \cite{H60})
of whether or not $\pc(2)=\frac12$, and the search for a rigorous proof
attracted a number of fine mathematicians into the field, including Lucio. Several important
partial results were proved, culminating in 1980 with Kesten's complete proof
that $\pc(2)=\frac12$, \cite{K82}.
The interest of the community then migrated towards the case $d \ge 3$, before returning
firmly to $d=2$ with the 2001 proof by Smirnov, \cite{SS01,SS01a}, of Cardy's formula.  

Lucio contributed a number of fundamental techniques to percolation theory 
during the period 1978--1988,
and the main purpose of the current paper is to describe these and to explore their significance.
We mention Russo's formula, the Russo--Seymour--Welsh (RSW)
inequalities, his study of
percolation surfaces in three dimensions, and of the uniqueness of the infinite open cluster, and finally 
Russo's approximate zero--one law. Russo's formula and RSW theory have proved of
especially lasting value in, for example, recent developments concerning
conformal invariance for critical percolation. 

In Section \ref{sec:isi}, we mention some of Lucio's results concerning percolation of $+/-$ spins
in the two-dimensional Ising model. It was quite a novelty in the 1970s to use percolation as a tool
to understand long-range order in the Ising model. Indeed, Lucio's work on the percolation model was 
motivated in part by his search for rigorous results in statistical mechanics. His approach 
to the Ising model has been valuable in two dimensions. 
In more general situations, the correct geometrical model has been recognised since to be
the \emph{random-cluster model} of Fortuin and Kasteleyn (see \cite{G-RCM}).

This short account is confined to Lucio's contributions to percolation, and does not touch on
his work lying closer to ergodic theory
and dynamical systems, namely \cite{CammR,LGR,FR81,MR}, and neither does
it refer to the paper \cite{FR}. A comprehensive list of Lucio's mathematical
publications, taken from MathSciNet, may be found at the end of this paper.

Results from Lucio's work
will be described here  using `modern' notation. No serious attempt is made to include
comprehensive citations of the related work of others.

\section{Russo's formula}\label{sec:form}

Let $\Om=\{0,1\}^E$ where $E$ is finite, and let $P_p$ be
product measure on the partially ordered set $\Om$ with density $p\in[0,1]$. 
An event $A\subseteq \Om$
is called \emph{increasing} if:
$$
\om\in A, \ \om\le\om' \q \Rightarrow\q  \om'\in A.
$$
Let $\om\in\Om$. An element $e\in E$ is called \emph{pivotal} for 
an increasing event $A$ if
$\om_e\notin A$ and $\om^e\in A$, where $\om^e$ and $\om_e$ are obtained
from $\om$ by varying the state of the edge $e$ thus:
$$
\om_e(f) = \begin{cases} 0 &\text{if } f=e,\\
\om(f) &\text{if } f \ne e,
\end{cases}\qq
\om^e(f) = \begin{cases} 1 &\text{if } f=e,\\
\om(f) &\text{if } f \ne e.
\end{cases}
$$
In other words, $e$ is said to be
pivotal for $A$ if the occurrence of $A$ depends on the state of $e$.
 
\begin{theorem}[Russo's formula, \cite{LR-crit}]
Let $A$ be an increasing event. We have that
$$
\frac d{dp} P_p(A) = \sum_{e\in E} P_p(e \text{\rm\ is pivotal for } A).
$$
\end{theorem}

Similar techniques are encountered independently
in related fields. For example, Russo's formula is essentially
equation (4.4) of Barlow and Proschan's book \cite[p.\ 212]{BaP} on reliability theory.
Such a formula appeared also in the work of Margulis, \cite{Marg}, in the Russian
literature. A characteristic of Lucio's work is the geometric context
of the formula when applied in situations such as percolation, and it is in this
context that Lucio's name is prominent. In a typical application to percolation,
one uses the geometrical characteristics of the event $\{e \text{  is pivotal for } A\}$ to
derive differential inequalities for $P_p(A)$. 

Russo's formula is key to the study of geometrical probability governed by 
a product measure. It has so many applications that is a challenge to single out any one.
We mention here its use in the derivation of exact values for critical exponents in two dimensions, \cite{Kes87,SW}. 

Similarly, extensions of Russo's formula have been central in several related fields, including
but not limited to the contact model \cite[Thm 2.13]{BG2},
continuum percolation \cite{FPR,JZG}, and the random-cluster model \cite[Prop.\ 4]{BGK}.

\section{Russo--Seymour--Welsh inequalities}\label{sec:RSW}

For twenty years from about 1960 to 1980, mathematicians attempted to prove that the critical
probability $\pc$ of bond percolation on the square lattice satisfies $\pc=\frac12$. This
prominent open problem was in the spirit of that of the critical temperature of 
the Ising model, resolved in 1944 by Onsager, \cite{Ons}. Harris \cite{H60} showed how
to use duality to obtain $\pc \ge \frac12$, but the corresponding upper
bound was elusive. Then, in 1978, a powerful technique  emerged in independent and contemporaneous
work of Lucio, \cite{LR-note}, and Seymour and Welsh, \cite{SeW}. It has come to
be known simply as `RSW'.

Consider bond percolation with density $p$ on the square lattice $\ZZ^2$.
A \emph{left--right crossing} 
of a rectangle $B$ is an open path
in $B$ which joins some vertex on its left side
to some vertex on its right side. For positive integers $m$ and $n$, we define
the rectangle
$$
B(m,n) = \bigl[0,2m\bigr]\times [0,2n],
$$
and let $\LR(m,n)$ be the event that there exists a left--right
crossing of $B(m,n)$. 

\begin{lemma}[Russo--Seymour--Welsh (RSW), \cite{LR-note,SeW}] \label{tcn70}
Let $p\in(0,1)$. We have that
$$
P_p\bigl( \LR(\tfrac32 n,n)\bigr)\geq \bigl(1-\sqrt{1-\tau}\bigr)^3,
$$
where $\tau=P_p(\LR(n,n))$.
\end{lemma}

This fundamental but superficially innocuous lemma implies that, 
if the chance of crossing a square is bounded from $0$ uniformly
in its size, then so is the chance of crossing a rectangle with aspect ratio $\frac32$.
Using the self-duality of $\ZZ^2$, we have as input to
the RSW lemma that 
\begin{equation}\label{eq:selfd}
P_{\frac12}(\LR(n,n))\ge \tfrac 12.
\end{equation}

Let $A_n$ be the event that the annulus $[-3n,3n]^2\setminus [-n,n]^2$
contains an open cycle with the origin in the bounded component of its complement in $\RR^2$.
Using elementary geometrical arguments and the FKG inequality, it follows by
the RSW lemma and \eqref{eq:selfd} 
that there exists $\sigma>0$ such that $P_p(A_n)\ge \sigma$ for $n \ge 1$
and $p\ge \frac12$.

The RSW lemma and the ensuing annulus inequality have proved to be key to the study
of percolation in two dimensions. In common with other useful methods of mathematics,
there is now a cluster of related inequalities, see for example \cite{BR}, \cite[Sect.\ 5.5]{G-pgs},
and \cite[Chap.\ 5]{Wer09}.

RSW methods were used by their discoverers to make
useful but incomplete progress towards proving that $\pc=\frac12$, and they played a role in Kesten's
full proof, \cite{K82}. (The principle novelty of Kesten's paper was a bespoke
theory of sharp threshold, see Section \ref{sec:01}.) 
More precisely, they led to the following result, which 
is presented in terms of \emph{site} percolation on the square lattice
$\ZZ^2$ and its \emph{matching lattice} 
$\ZZ^2_\ast$, derived by adding the two diagonals to each face of $\ZZ^2$.

 \begin{theorem}[Russo, \cite{LR-note}]\label{thm:rus2}
 Consider site percolation on the square lattice $\ZZ^2$. The critical points
 $$
 \pc=\sup\{p: \theta(p)=0\},\qq \pic=\sup\{p: \chi(p)<\oo\},
 $$
 satisfy
 \begin{equation}\label{eq:pcrit}
 \pc+\pic^*=1,\qq \pc^*+\pic=1,
 \end{equation}
 where an asterisk denotes the corresponding values on the matching lattice.
 \end{theorem} 
 
The parallel work of Seymour and Welsh, \cite{SeW}, was directed at the bond model on $\ZZ^2$, of which
the dual model lies on a translate of $\ZZ^2$.  
Following Kesten's proof of $\pc(\ZZ^2)=\frac12$ for bond percolation, 
Lucio revisited Theorem \ref{thm:rus2} in \cite{LR-crit} with a proof that 
$\pic^*=\pc^*$, and the consequent improvement of \eqref{eq:pcrit},
namely $\pc+\pc^*=1$. He also completed the proof, begun in \cite{LR-note},
that $\theta$ (and, similarly,
the dual percolation probability $\theta^*$) is a continuous function
on $[0,1]$. Continuity in two dimensions has since been extended to
general percolation models (see, for example, \cite[Sect.\ 8.3]{G99}). 

RSW theory is now recognised as fundamental to rigorous proofs of
conformal invariance of critical two-dimensional percolation and all that comes with that. 
The proof of Cardy's formula, \cite{SS01,SS01a}, provides a major illustration. 
It was observed by Aizenman and Burchard, \cite{AizB}, that certain connection probabilities
belong to a space of uniformly H\"older functions. Since this space is compact, such
functions have subsequential limits as the mesh of the lattice approaches $0$.
The above H\"older property is proved using annulus inequalities.
 
Indeed the power of RSW arguments extends beyond percolation to a  host of problems
involving two-dimensional stochastic geometry, such as the FK-Ising model \cite{CDCH}
and Voronoi percolation \cite{Tass15}. In addition, RSW theory provides one of the main
techniques for the proof by Beffara and Duminil-Copin, \cite{BDC}, that
the random-cluster model on $\ZZ^2$ with cluster-weighting parameter $q\ge 1$
has critical value $\pc(q)=\sqrt q/(1+\sqrt q)$. We retrieve Kesten's theorem by setting $q=1$.

\section{Approximate zero--one law}\label{sec:01}

Kolmogorov's zero--one law may be stated as follows.
Consider the infinite product space $\Om=\{0,1\}^\NN$ endowed
with the product $\sigma$-algebra and the product measure $P_p$.
If $A$ is an event that is independent of any finite subcollection
$\{\om(e): e \in E\}$, $E \subseteq \NN$, $|E|<\oo$, then
$P_p(A)$ equals either $0$ or $1$. It follows that, for
an increasing event $A$,  there exists $p_0\in[0,1]$ such that
\begin{equation*}
P_p(A) = \begin{cases} 0 &\text{if }  p<p_0,\\
1 &\text{if } p>p_0.
\end{cases}
\end{equation*}
This law is intrinsically an \emph{infinite-volume} effect, in that the index set is the
infinite set $\NN$.
Lucio posed the farsighted question in \cite{LR-app} of whether there exists
a \emph{finite volume} version of this result, and this led him
to his `approximate zero--one law', following.

Let $\Om=\{0,1\}^E$ where $E$ is finite, and let $P_p$ be
product measure on the partially ordered set $\Om$ with density $p\in[0,1]$. 
The \emph{influence} $I_{A,p}(e)$ of $e\in E$ on the event $A\subseteq \Om$
is defined by
$$
I_{A,p}(e) =P_p\bigl(1_A(\om_e)\ne 1_A(\om^e)\bigr),
$$
where $1_A$ denotes the indicator function of $A$. When $A$ is
increasing, this may be written
\begin{equation}\label{-1}
I_{A,p}(e) =P_p(\om_e \notin A, \ \om^e\in A)=
P_p(e \text{ is pivotal for } A).
\end{equation}

\begin{theorem}[Russo's approximate zero--one law, \cite{LR-app}]\label{thm:app}
For  $\eps>0$, there exists $\eta>0$ such that, if $A$ is an increasing
event and
\begin{equation}\label{eq:as}
 I_{A,p}(e)<\eta,\qq  e \in E,\ p\in [0, 1],
\end{equation}
then there exists $p_0\in[0,1]$ such that
\begin{equation}\label{eq:conc}
P_p(A) \begin{cases} \le\eps &\text{if }  p<p_0-\eps,\\
\ge 1-\eps &\text{if } p>p_0+\eps.
\end{cases}
\end{equation}
\end{theorem}

This result was motivated by a desire to generalise certain results for box-crossing probabilities
in percolation.
Its impact extends far beyond percolation, and it is a precursor of a more recent
theory, pioneered by Kahn, Kalai, Linial, \cite{KKL} and Talagrand, \cite {Tal94,Tal99}, of influence
and sharp threshold.
It is proved at \cite[Thm 1.1]{Tal94} that there exists
an absolute constant $c>0$ such that, for $p\in(0,1)$ and an increasing event $A$,
\begin{equation}\label{0}
\sum_{e\in E} I_{A,p}(e) \ge \left(\frac{c}{p(1-p)
\log[2/(p(1-p))]}\right)P_p(A)(1-P_p(A))\log(1/m_p),
\end{equation}
where
$$
m_p=\max\{I_{A,p}(e): e \in E\}.
$$
It follows that, when $p\in(0,1)$ and $A$ is increasing,
\begin{equation}\label{000}
\sum_{e\in E} I_{A,p}(e) \ge c' P_p(A)(1-P_p(A))\log(1/m_p),
\end{equation}
where $c'>0$ is an absolute constant.

Amongst the  implications of \eqref{000} is a quantification of the relationship
between $\eps$ and $\eta$ in Theorem \ref{thm:app}.
Suppose \eqref{eq:as} holds with $\eta\in(0,1)$, 
so that $m_p\le \eta$. By \eqref{000} and Russo's formula,
\begin{equation}\label{eq:int}
\frac d{dp}P_p(A)  \ge  c'P_p(A)(1-P_p(A))\log (1/\eta).
\end{equation}
Choose $p_0$ such that $P_{p_0}(A)=\frac12$, and integrate
\eqref{eq:int} to obtain 
\begin{equation}\label{eq:conc1}
P_p(A) \begin{cases} \le\eps_p &\text{if }  p<p_0,\\
\ge 1-\eps_p &\text{if } p>p_0,
\end{cases}
\end{equation}
with
$$
\eps_p=\frac1{1+(1/\eta)^{c'|p-p_0|}}.
$$

Such inequalities have found numerous applications in percolation and related topics, 
see for example \cite{BR,DCM,GrGr2}. 
They have been extended to general product measures, \cite{BKKKL,GJN},
and to probability measures satisfying the FKG lattice condition, \cite{GrGr1}.
Recent overviews include \cite[Chap.\ 4]{G-pgs} and \cite{KalS}.

\section{Percolation in dimension $d\ge 3$}\label{sec:d3}

In 1983, Lucio spent a sabbatical at Princeton University. His work
during that period led to two significant publications \cite{ACCFR,CR} 
on aspects of percolation in three dimensions. The first of these
caused quite a stir in the community at the time of its appearance, largely
since most work until then had been for models in only two dimensions.
Whereas the dual of a bond model in two dimensions is another bond model,
the dual model in three or more dimensions is a `plaquette' model. Since the topology of
surfaces of plaquettes is much more complicated that that of paths, the ensuing
percolation duality poses a number of challenging topological questions.  

The authors of \cite{ACCFR} consider bond percolation on $\ZZ^3$
with density $p$, together with its dual `plaquette'
model on $\ZZ^3_*:=\ZZ^3+(\frac12,\frac12,\frac12)$.
A \emph{plaquette} is a unit square with vertices in $\ZZ^3_*$, and   its bounding lines
are \emph{edges} of $\ZZ^3_*$. Each edge $e$ of $\ZZ^3$ 
intersects a unique plaquette $\Pi_e$, and $\Pi_e$ is termed
\emph{occupied} if and only if $e$ is closed (and \emph{unoccupied}
otherwise). Thus, a plaquette is occupied
with probability $1-p$. For any  collection $F$ of plaquettes, the boundary
$\pd F$ is defined to be set of edges of $\ZZ^3_*$  belonging to an \emph{odd}
number of members of $F$.

Let $\gamma$ be a cycle of $\ZZ^3_*$.
The main results of \cite{ACCFR} concern the probability there exists a
set $F$ of occupied plaquettes which spans $\gamma$ in the sense that $\pd F = \gamma$. 
For simplicity, we shall suppose here that $\gamma$ is a $m \times n$ rectangle
of the $x/y$ plane, and we denote the above event as $W_\gamma$.
Note that $\gamma$ has area $mn$ and perimeter $2(m+n)$.

\begin{theorem}[Aizenman, Chayes, Chayes, Fr\"ohlich, Russo, \cite{ACCFR}]
There exist constants $\pic,\rhoc\in(0,1)$ such that
\begin{equation*}
-\log P_p(W_\gamma) \sim \begin{cases}
\alpha mn &\text{if } 1-p <\pic,\\
\beta (m+n) &\text{if } 1-p>\rhoc,
\end{cases}
\end{equation*}
where $\alpha, \beta>0$ depend on $p$, and the asympotic relation is as $m,n\to\oo$.
\end{theorem}

The constants $\pic$, $\rhoc$ are the critical densities  of the bond percolation model
on $\ZZ^3$ given by
\begin{align*}
\pic= \sup\{p:\chi(p)<\oo\},\qq \rhoc=\lim_{k\to\oo}\wpc(k),
\end{align*}
where $\wpc(k)$ is the \emph{slab critical point}
$$
\wpc(k) = \sup\bigl\{p: P_p\bigl(0\lra \oo \text{  in } [0,\oo)^2\times[0,k]\bigr)=0\bigr\}.
$$
It was conjectured in \cite{ACCFR} that $\pic=\pc=\rhoc$. The first equality was proved
later in \cite{AB,Men}, and the second in \cite{BGN,GrM}.

There are only few percolation models on finite-dimensional lattices
for which the numerical values of the critical probabilities are known 
exactly, and all such exact results are in two dimensions only (see, for example, \cite{GM14}).
In contrast, quite a lot of work
has been devoted to obtaining rigorous upper and lower bounds
for critical probabilities, and there is a host of numerical
estimates.

Consider site percolation on the simple cubic lattice $\ZZ^3$.
By a comparison with the site model on the triangular lattice, Lucio has shown (with 
Campanino, in \cite{CR}) that
$\pc\le \frac12$. (See also \cite{Men85}.) 
They obtained also the strict inequality, with a distinctly more complicated argument. 

\begin{theorem}[Campanino, Russo, \cite{CR}]\label{thm:camR}
The critical probability of site percolation on $\ZZ^3$ satisfies $\pc < \frac12$.
\end{theorem}

The point of this work was to show that, in a neighbourhood of $p=\frac12$,
there is coexistence of infinite open and infinite closed clusters in $\ZZ^3$.
The corresponding statement for $d=2$ is, of course, false, in that coexistence occurs for
no value of $p$.

Theorem \ref{thm:camR} may still be the best rigorous
upper bound that is currently known for $\pc$. By examining its proof, 
one may calculate a small $\eps>0$ such that $\pc<\frac12-\eps$. It is expected that 
$\pc \approx 0.31$.

\section{Uniqueness of the infinite open cluster}\label{sec:unique}

Let $I$ be the number of infinite open clusters of a percolation model
in a finite-dimensional space. For a period
in the 1980s, the `next' problem was to prove that $P_p(I=1)=1$ in the supercritical
phase (when $p>\pc$). This problem was solved by
Aizenman, Kesten, and Newman \cite{AKN} in 1987. Their proof
seemed slightly mysterious at the time, and it was simplified by Lucio in the jointly
written paper \cite{GGR}. The key step was to show, using a  large-deviation estimate
present already in \cite{AKN}, that there is density $0$ of sites that are adjacent to \emph{two} 
distinct infinite clusters. 

This useful argument was soon overshadowed by the beautiful proof of uniqueness by Burton 
and Keane, \cite{BK}, of which a key step is a novel argument to show
there is density $0$ of sites that are adjacent to \emph{three} 
distinct infinite clusters. The proof of \cite{BK} uses translation-invariance of the underlying measure
together with a property of so called `finite energy',
and may thus be extended to more general measures than product measures.
On the other hand, since the proof uses no quantitative estimate, it yields no `rate'.
The methods of \cite{AKN,GGR} provide a missing rate, and this has been useful in the 
later work \cite{Cerf15,ChattS}. 

The question of uniqueness for dependent models is potentially harder, 
since the large-deviation estimate of \cite{AKN,GGR} is not available.
In joint work \cite{GKR} with Gandolfi and Keane, Lucio used path-intersection
arguments to show uniqueness for ergodic, positively associated measures in two dimensions,
satisfying certain translation and reflection symmetries.  Unlike
the Burton--Keane proof, they needed no finite-energy assumption.
An application of this work to quantum spin systems may be found in \cite{AN}.

\section{Ising model}\label{sec:isi}

Lucio has written three papers on the geometry of the $d$-dimensional
Ising model, \cite{CNPR2,CNPR1,LR-inf}. In this work, he (and his coauthors)
studied the relationship between properties of the infinite-volume Gibbs measures
and the existence or not of an infinite cluster of either $+$ or $-$ spins (that is,
of percolation in the Ising model).

The first two of these papers \cite{CNPR2,CNPR1} explore a relationship between the 
Ising magnetization and the above percolation probability, and yield the
non-existence of percolation
in the high-temperature phase.
This is complemented when $d=2$ with the proof that percolation (of the
corresponding spin) exists in 
the low-temperature phase for either of the pure infinite-volume limits $\mu_+$, $\mu_-$, 
obtained \resp\ as the weak limits with $+/-$ boundary conditions.
These methods were developed further in \cite{CNPR1}, where a phase diagram was proposed for the 
existence of infinite clusters in the two-dimensional ferromagnetic Ising model,
as a function of external field $h$ and temperature $T$. The principal features
of this diagram were later proved  by Higuchi, \cite{Hig93}.

One of the central problem in two dimensions of the late 1970s was to prove or disprove the
statement that every infinite-volume Gibbs measure is a convex combination of the 
two extremal measures $\mu_+$, $\mu_-$.
Lucio obtained the following important result for this problem. 

\begin{theorem}[Russo, \cite{LR-inf}]\label{thm:convex}
Any infinite-volume Gibbs measure $\mu$, which is trans\-lation-invariant in one or both 
of the axial directions, is a convex combination of $\mu_+$ and $\mu_-$.
\end{theorem}

Lucio proved this by considering the existence (or not)
of infinite $+/-$ clusters on $\ZZ^2$ and its matching lattice.
The full conclusion, without an assumption of partial translation-invariance, was obtained later 
in independent work of  Aizenman, \cite{Aiz80}, and Higuchi, \cite{Hig81} (see also \cite{GHig}).
Therefore, in two dimensions (unlike three dimensions)
there exists no non-translation-invariant Gibbs measure.

More recent work on the geometrical properties of the Ising model has been
centred around the random-cluster model and the random-current representation, rather than 
the more fundamental percolation model. See, for example, \cite{ADCS,G-RCM}.

\section*{Acknowledgement} 
This work was supported in part by the Engineering
and Physical Sciences Research Council under grant EP/I03372X/1.
The author thanks Alberto Gandolfi for his comments and suggestions.

\renewcommand\refname{Mathematical publications of Lucio Russo}

\bibliographystyle{amsplain}
\providecommand{\bysame}{\leavevmode\hbox to3em{\hrulefill}\thinspace}
\providecommand{\MR}{\relax\ifhmode\unskip\space\fi MR }
\providecommand{\MRhref}[2]{%
  \href{http://www.ams.org/mathscinet-getitem?mr=#1}{#2}
}
\providecommand{\href}[2]{#2}

\renewcommand\refname{References}
\providecommand{\bysame}{\leavevmode\hbox to3em{\hrulefill}\thinspace}
\providecommand{\MR}{\relax\ifhmode\unskip\space\fi MR }
\providecommand{\MRhref}[2]{%
  \href{http://www.ams.org/mathscinet-getitem?mr=#1}{#2}
}
\providecommand{\href}[2]{#2}


\end{document}